\documentclass[preprints,article,submit,moreauthors,pdftex,10pt,a4paper]{mdpi}
\firstpage{1}
\pubvolume{11}
\issuenum{1}
\articlenumber{5}
\pubyear{2019}
\copyrightyear{2019}
\history{Received: 5 February 2019; Accepted: 12 March 2019; Published: date}
\usepackage{amsmath}
\usepackage{color,soul}
\usepackage{diagbox}
\usepackage{booktabs}
\usepackage{titletoc}
\renewcommand\thesubsubsection{\Alph{subsubsection}}

\setitemize{parsep=6pt,itemsep=0pt,leftmargin=*,labelsep=5.5mm}
\setenumerate{parsep=6pt,itemsep=0pt,leftmargin=*,labelsep=5.5mm}
\setlist[description]{itemsep=0mm}

\Title{A Semi-Markov Model for the Control of Thermostatically Controlled Load}

\Author{Benyuan Zhao, Peichao Zhang and Yizhi Cheng}

\AuthorNames{Benyuan Zhao and Peichao Zhang and Yizhi Cheng}

\address[1]{%
Department~of~Electrical~Engineering, Shanghai Jiao Tong University, Shanghai 200240, China; origin\_zhao@163.com (B.Z.); tabdel@sjtu.edu.cn (Y.C.)}

\corres{\hangafter=1 \hangindent=1.05em \hspace{-0.82em}Correspondence: pczhang@sjtu.edu.cn}

\abstract{Considering the potential of thermostatically controlled loads (TCLs) to provide flexibility in demand response or load control, a semi-Markov model (SMM) for the ON/OFF controlled TCL is developed in this paper. This model makes full use of the adjustment flexibility of TCLs when the control period is long and maintains the diversity of switch state in the cluster. This model also can satisfy user comfort and protect user privacy. Then, this paper adopts the cyber-physical system (CPS) to realize the coupling of the discrete control process and the continuous physical process. Finally, the proposed model is applied to the coordination of large-scale heterogenous air-conditioners (ACs) based on the equivalent thermal parameters (ETP) model. Simulation results verify that under the proposed approach, the power of TCLs cluster can track the control signal accurately, with both user comfort and diversity of TCL cluster's operation states guaranteed.}

\keyword{Thermostatically controlled load (TCL), semi-Markov control model (SMM), cyber-physical system (CPS), equivalent thermal parameters (ETP)}

\begin{document}

\section{Introduction}
The ever increasing penetration of renewable energy sources (RES) in the power grid has posed a great threat to the power system on meeting the supply-demand balance. Considering the intermittency and uncertainty brought by the integration of RES, traditional generators of slow response speed and high regulation cost may fail to ensure the reliability and efficiency of the power system. At the same time, the consumption of thermostatically controlled loads (TCLs), i.e., air-conditioners (ACs), refrigerators and water heaters, has grown rapidly in China recently. Especially during summer, ACs would consume about $30{\rm \%} \sim 40{\rm \%}$, and even exceed ${50 \rm\% }$ of total electricity in some large and medium-sized cities. Due to their potential to provide flexibility by changing the operating cycle and adjusting the temperature setting value, TCLs can participate in various ancillary services while ensuring user comforts, such as active power balance, frequency regulation, and peak shaving.

To fully explore the flexibility of TCLs, it is of great significance to accurately model their thermal dynamics. At present, a large proportion of TCLs are under ON/OFF control that their instantaneous power cannot be adjusted continuously, with the commonly adopted control strategies being switching control, temperature control, and probability control. Switching control refers to the strategy that controls an aggregation of TCLs to track a target power by directly controlling the switching state of each TCL. In \cite{lu2004state} a state-queueing model is used to divide the TCLs into clusters according to their indoor temperature and switching state, with the number of TCLs to be turned on determined by target power. In \cite{hao2014aggregate}, the state of TCL is controlled based on the relationship between the target power and the power of TCL in a free state. Temperature control refers to the control of TCL by adjusting its set point. In \cite{WANGChengshan2012,WANGDan2014}, the temperature control method is combined with the state sequence model of the cluster to achieve accurate tracking of the control signal. However, since both switching control and temperature control are direct control strategies, the control privileges of TCLs should be acquired by controllers, which may cause problems concerning information security and potential load oscillation.

To address these problems, a probability control that strategically controls the duty cycle of TCL to track the average power of a cluster of TCLs are proposed in many studies with the adoption of the Markov model. In \cite{paccagnan2015range}, the relationship between the probability distribution of cluster state and transition probability is established by the Fokker-Planck equation to realize the control of a TCL cluster. To facilitate control, the temperature distributions in the cluster are divided into several groups in Literature \cite{mathieu2012state,zhang2012aggregate}. Literature \cite{totu2016demand} uses the Markov chain to describe the dynamic process of TCL's number change in each interval. Literature \cite{moya2014hierarchical,williams2016control,kalsi2015loads} proposes a two-layer controlled method. Such a method obtains duty cycle of the cluster by frequency adjustment on the upper layer and computes transition probability through the change of the duty cycle. For the multiple solutions of transition probability, the switching times have been optimized for the only solution in Literature \cite{williams2016control,kalsi2015loads}.

 As a kind of indirect load control method, the probability control method realized by the local controller uses the transition probability as a control signal and increases the clusters' randomness to ensure the diversity of TCLs' state. However, the aforementioned approaches based on the Markov chain have several limitations that need to be addressed. First, the aggregation process proposed in \cite{paccagnan2015range} requires a collection of model parameters and user preferences that may violate users' privacy. Second, according to the aggregation methods proposed in \cite{paccagnan2015range,mathieu2012state}, they can only be applied to the coordination of TCLs with homogeneous parameters instead of heterogenous TCLs. Besides, methods that omit the lock time constraint \cite{mathieu2012state,zhang2012aggregate} may arise frequent switches of the TCLs, thus increasing wear and causing damage to equipment. Although this constraint is considered in \cite{moya2014hierarchical,williams2016control,kalsi2015loads}, these methods are not suitable for applications with the control period longer than the lock time, to fully explore the flexibility that TCLs can offer. 

This paper proposes a probability controlled model based on the semi-Markov process that can systematically address the aforementioned problems. Besides, the cyber-physical system is adopted to realize the autonomous control of TCL by establishing the coupling relationship between the physical process and the control process.

The rest of this paper is organized as follows. Section \ref{sec:Semi-Markov Control Model} introduces the probability control model for ON/OFF controlled TCLs based on the semi-Markov model. In Section \ref{sec:Physical Model}, based on the theory of CPS, the SMM is applied to control a population of ACs whose physical model is derived from the ETP model. Section \ref{sec:case_studies} develops case studies to prove the effectiveness of the method. Finally, some concluding remarks are given in Section \ref{sec:conclusions}.

\section{Control Model}
\label{sec:Semi-Markov Control Model}

\subsection{Semi-Markov Model}
\label{sec:SMM}
\begin{figure}[H]
 \centering
 \includegraphics[width=8 cm]{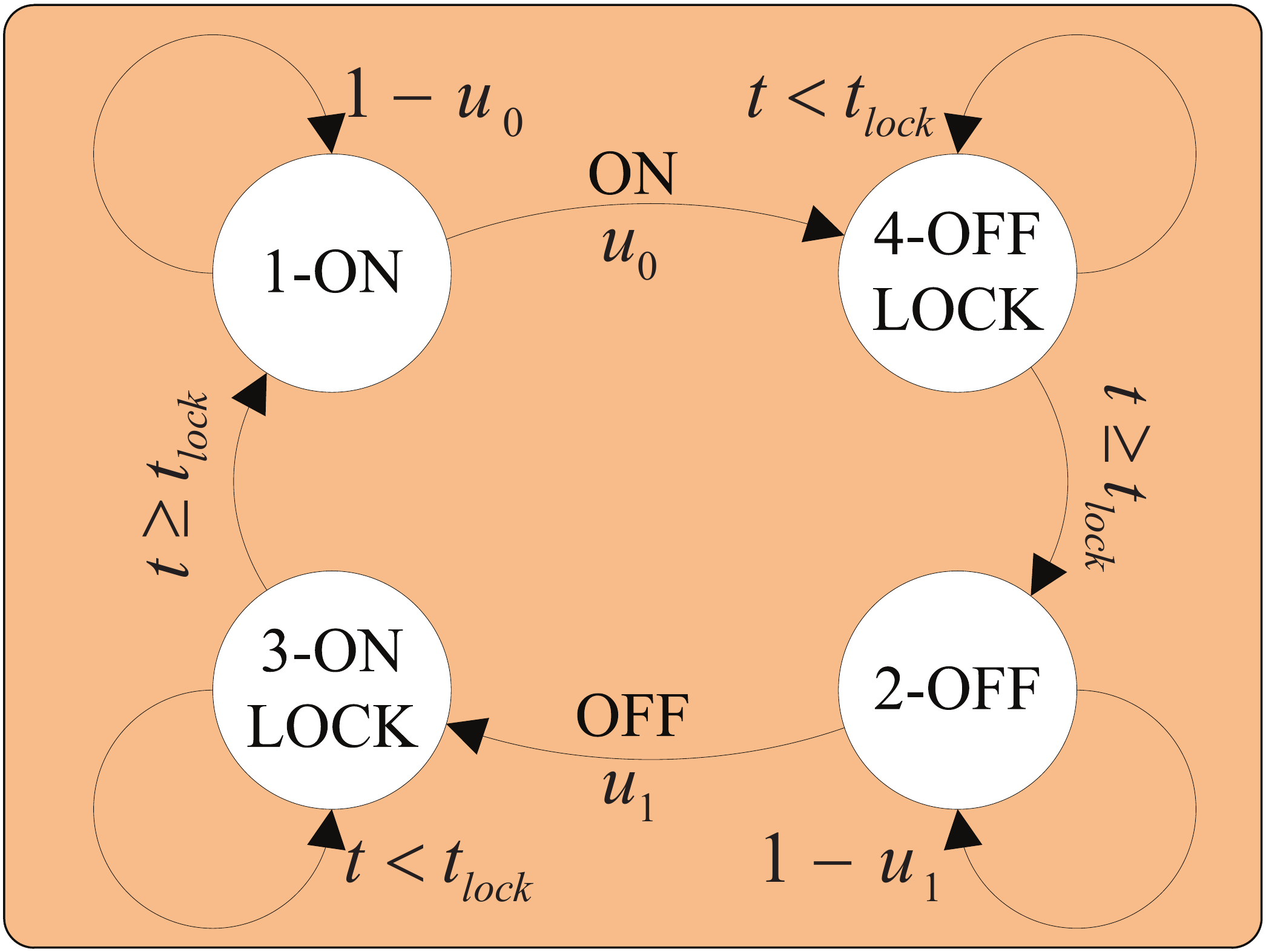}
 \caption{Semi-Markov model of the ON/OFF controlled TCL.}
 \label{fig:semi-markov}
\end{figure}

The state machine of an ON/OFF controlled TCL is demonstrated in Figure \ref{fig:semi-markov}. As shown, four states are included after the introduction of lock time to prevent the frequent start and stop. Both $ON$ and $ONLOCK$ are the Power-on state of TCL, and the difference is that the former satisfies the lock time and can be OFF at any time; both $OFF$ and $OFFLOCK$ are the Power-off state of TCL, and the difference is that the former satisfies the lock time and can be ON at any time. 

The Markov process requires that the time distribution of the state transition of the system be memoryless. However, in Figure \ref{fig:semi-markov}, the transition duration of the lock state does not obey the exponential distribution due to the constant lock time. To this end, this paper proposes a semi-Markov model (SMM) to establish the generalized control model for all ON/OFF controlled TCLs that uses two transition probabilities $u_0$ and $u_1$ to control the transition between states. 

The SMM is defined as follows \cite{harlamov2008continuous}:
\begin{enumerate}
 \item [(1)] In the state space $\aleph {\rm{ = \{ }}1,2,3,4\}$, the probability that the system leave state $m$ to state $n$ is ${{p_{mn}}}$, thus
 \begin{equation}
 \label{eq:probability_m_n}
 \sum\limits_{n \ne m,n = 1}^r {{p_{mn}} = 1}\,.
 \end{equation}
 \item [(2)] Given that the next state is $n$, the transition time from $m$ to $n$ satisfies the distribution ${{F_{mn}}(t)}$, which can be any distribution. So the distribution ${F_m}(t)$ in state $m$ is
 \begin{equation}
 \label{eq:distribution_Fmn}
 {F_m}(t) = \sum\limits_{m \ne n,m = 1}^r {{p_{mn}}{F_{mn}}(t)}\,.
 \end{equation}

 The average stay time of state $m$ is
 \begin{equation}
 \label{eq:stay_time_of_m}
 {T_m} = \int_0^\infty {td{F_m}(t)}\,.
 \end{equation}

 Let $m = 1, 2, 3, 4$ represent $ON$, $OFF$, $ONLOCK$ and $OFFLOCK$, respectively. Figure \ref{fig:semi-markov} shows that the transfer of state is
 unidirectional, so according to Equation \eqref{eq:probability_m_n},
 \begin{equation}
 \label{eq:state_transfer_of_FFA}
 {p_{14}} = {p_{42}} = {p_{23}} = {p_{31}} = 1\,.
 \end{equation}

 ${F_{14}}(t)$ and ${F_{23}}(t)$ are exponential distributions with parameters ${\frac{{{u_0}}}{{\Delta t}}}$ and ${\frac{{{u_1}}}{{\Delta t}}}$ ,respectively. Thus
 \begin{equation}
 \label{eq:F14_and_F23}
 {F_{14}}(t) = 1 - {e^{ - \frac{{{u_0}}}{{\Delta t}}t}},{F_{23}}(t) = 1 - {e^{ - \frac{{{u_1}}}{{\Delta t}}t}}\,,
 \end{equation}
 where ${\Delta t}$ is the interval of semi-Markov process.

 According to Equation \eqref{eq:stay_time_of_m}, the average stay of state 1, 2 and its standard deviation are
 \begin{equation}
 \label{eq:average_stay_time_of_state12}
 \left\{ \begin{array}{l}
 {T_1}{\rm{ = }}\int_0^\infty {\frac{{{u_0}}}{{\Delta t}}t{e^{ - \frac{{{u_0}}}{{\Delta t}}t}}dt = } \frac{{\Delta t}}{{{u_0}}}\\
 {T_2}{\rm{ = }}\int_0^\infty {\frac{{{u_1}}}{{\Delta t}}t{e^{ - \frac{{{u_1}}}{{\Delta t}}t}}dt = } \frac{{\Delta t}}{{{u_1}}}
 \end{array} \right.\,.
 \end{equation}

 \begin{equation}
 \label{eq:standard_deviation}
 \left\{ \begin{array}{l}
 {\sigma _1}{\rm{ = }}\int_0^\infty {\frac{{{u_0}}}{{\Delta t}}{t^2}{e^{ - \frac{{{u_0}}}{{\Delta t}}t}}dt} - {\left( {\int_0^\infty {\frac{{{u_0}}}{{\Delta t}}t{e^{ - \frac{{{u_0}}}{{\Delta t}}t}}dt} } \right)^2} = \frac{{\Delta t}}{{{u_0}}}\\
 {\sigma _2}{\rm{ = }}\int_0^\infty {\frac{{{u_1}}}{{\Delta t}}{t^2}{e^{ - \frac{{{u_1}}}{{\Delta t}}t}}dt} - {\left( {\int_0^\infty {\frac{{{u_1}}}{{\Delta t}}t{e^{ - \frac{{{u_1}}}{{\Delta t}}t}}dt} } \right)^2} = \frac{{\Delta t}}{{{u_1}}}
 \end{array} \right.\,.
 \end{equation}

 Obviously, the average stay of state3, 4 are
 \begin{equation}
 \label{eq:average_stay_time_of_state34}
 {T_3} = {T_4} = {t_{{\rm{lock}}}}\,.
 \end{equation}
 \item [(3)]If the semi-Markov chain is irreducible, after a finite time, the probability of being in a certain state can converge to a constant independent of the initial state , denoted as steady-state probability \cite{harlamov2008continuous}, while the probability before convergence is called transient probability. According to the characteristic of the semi-Markov chain, the steady-state probability ${p_m}$ is
 \begin{equation}
 \label{eq:probability_of_stability_Pm}
 {p_m} = \frac{{{\pi _m}{T_m}}}{{\sum\limits_{m = 1}^r {{\pi _m}{T_m}} }}\,,
 \end{equation}
 where ${{\pi _m}}$ represents the distribution of state $m$ that satisfies the equation
 \begin{equation}
 \label{eq:stationary_of_state_m}
 \left\{ \begin{array}{l}
 {\pi _n} = \sum\limits_{m = 1}^r {{\pi _m}{p_{mn}}} \\
 \sum\limits_{m = 1}^r {{\pi _m} = 1}
 \end{array} \right.\,.
 \end{equation}

 Substitute Equation \eqref{eq:state_transfer_of_FFA} into Equation \eqref{eq:stationary_of_state_m}, it can be obtained that
 \begin{equation}
 \label{eq:pi}
 {\pi _1} = {\pi _2} = {\pi _3} = {\pi _4} = {\rm{0}}{\rm{.25}}\,.
 \end{equation}

 Combine Equation \eqref{eq:probability_of_stability_Pm} with Equation \eqref{eq:pi}, the steady-state probability is
 \begin{equation}
 \label{eq:probability_of_stability}
 {p_m} = \frac{{{T_m}}}{{{T_1} + {T_2} + {T_3} + {T_4}}}{\kern 1pt} {\kern 1pt} {\kern 1pt} {\kern 1pt} {\kern 1pt} {\kern 1pt} {\kern 1pt} {\kern 1pt} {\kern 1pt} {\kern 1pt} {\kern 1pt} {\kern 1pt} {\kern 1pt} {\kern 1pt} m = 1,2,3,4\,.
 \end{equation}

 It can be inferred that the steady-state probability is determined by the transition probabilities ${{u_0}}$ and ${{u_1}}$.

\end{enumerate}
\subsection{Determination of Transition Probabilities}
 For TCL $i$, its instantaneous power $\tilde P_{i}^k{\rm}$ is a discrete random variable, i.e., either $P_{i}^{rate}$ which represents the rated power or $0$, and its expected power during a control period $k$ denoted as $P_{\mathrm{exp}.i}^k$, is calculated as
\begin{equation}
\label{eq:expected_value_of_PFFA}
P_{\mathrm{exp}.i}^k = E(\tilde P_{i}^k{\rm{) = (}}{p_{\rm{1}}}{\rm{ + }}{p_{\rm{3}}}{\rm{)}}P_{i}^{rate}\,.
\end{equation}
Therefore, the expected power $P_{\mathrm{exp}.i}^k$ is also determined by the transition probabilities.

Assuming that the target power received at control period $k$ is $P_{i}^k$, let
\begin{equation}
\label{eq:expected_value_of_PFFA_equal_response_power}
P_{\mathrm{exp}.i}^k = P_{i}^k\,.
\end{equation}
According to Equation \eqref{eq:expected_value_of_PFFA}\eqref{eq:expected_value_of_PFFA_equal_response_power}, the solution of ${{u_0}}$ and ${{u_1}}$ is not unique. Considering the convergence performance during transient phase, their value should be determined carefully.

On one hand, since the actual durations of state ${\rm{1-ON,2-OFF}}$ would fluctuate seriously around the expected value within one cycle, TCLs with the indoor temperature near the boundary are prone to exceed the comfort limit. Therefore, according to Equation \eqref{eq:standard_deviation}, ${{u_0}}$ and ${{u_1}}$ should be as large as possible to decrease the standard deviation of duration time. Considering that given a constant ${{p _m}}$, ${{u_0}}$ and ${{u_1}}$ would increase or decrease simultaneously according to Equation \eqref{eq:probability_of_stability}, the following constraint are added:
\begin{equation}
\label{eq:the_added_value}
\left\{ \begin{array}{l}
{u_{\rm{1}}} = 1{\kern 1pt}, \;\;\mathrm{if} \;{{P_{i}^k} \mathord{\left/
 {\vphantom {{P_{i}^k} {P_{i}^{rate}}}} \right.
 \kern-\nulldelimiterspace} {P_{i}^{rate}}} > 0.5\\
{u_{\rm{0}}} = 1{\kern 1pt}, \;\;\mathrm{if} \;{{P_{i}^k} \mathord{\left/
 {\vphantom {{P_{AC.i}^k} {P_{i}^{rate}}}} \right.
 \kern-\nulldelimiterspace} {P_{i}^{rate}}} \le 0.5
\end{array} \right.\,.
\end{equation}

On the other hand, according to the convergence of the semi-Markov process, the state space cannot be a simple loop indicating that the transition loop cannot be $A$ to $B$ and then $B$ to $A$. When ${{P_{i}^k} \mathord{\left/
 {\vphantom {{P_{i}^k} {P_{i}^{rate}}}} \right.
 \kern-\nulldelimiterspace} {P_{i}^{rate}}} $ in Equation \eqref{eq:the_added_value} gradually reduces to 0.5, the expected duration time $T_2=\Delta t$ and ${{T_1}}$ decreases until it is finally negligible compared to ${{t_{{\rm{lock}}}}}$. On such condition, each TCL switches between $3-ONLOCK$ and $4-OFFLOCK$ every ${{t_{{\rm{lock}}}}}$ and the semi-Markov process is a simple loop. So in order to prevent this phenomenon, the constraint in Equation \eqref{eq:the_added_value} is replaced with a lower limit for $T_1$
 \begin{equation}
 \label{eq:T1limit}
 {T_1} > t\,,
 \end{equation}
 where $t$ is set to $1min$ in this paper, given that the lock time is typically $3-5min$ \cite{zhang2013aggregated}.

 It can proved that Equation\eqref{eq:T1limit} holds when the following condition satisfies
\begin{equation}
\label{eq:T2_fixed}
{{P_{i}^k} \mathord{\left/
 {\vphantom {{P_{i}^k} {P_{i}^{rate}}}} \right.
 \kern-\nulldelimiterspace} {P_{i}^{rate}}} > \frac{{t{\rm{ + }}{t_{{\rm{lock}}}}}}{{t + {t_{{\rm{lock}}}}{\rm{ + }}{\Delta t}{\rm{ + }}{t_{{\rm{lock}}}}}} > 0.5\,.
\end{equation}

Similarly, when ${{P_{i}^k} \mathord{\left/
 {\vphantom {{P_{i}^k} {P_{i}^{rate}}}} \right.
 \kern-\nulldelimiterspace} {P_{i}^{rate}}} $ in Equation \eqref{eq:the_added_value} gradually increases to 0.5, $T_2 > t$ would hold as long as the following condition is satisfied
\begin{equation}
\label{eq:T1_fixed}
{{P_{i}^k} \mathord{\left/
 {\vphantom {{P_{i}^k} {P_{i}^{rate}}}} \right.
 \kern-\nulldelimiterspace} {P_{i}^{rate}}} < \frac{{\Delta t + {t_{{\rm{lock}}}}}}{{\Delta t + {t_{{\rm{lock}}}}{\rm{ + }}t+{t_{{\rm{lock}}}}}} < 0.5\,.
\end{equation}

On the contrary, when neither \eqref{eq:T1_fixed} nor \eqref{eq:T2_fixed} is satisfied, another constraint is added to \eqref{eq:the_added_value} as follows to guarantee the lower limit for duration time
\begin{equation}
\label{eq:value_of_u0_and_u1}
\left\{ \begin{array}{l}
{u_1} = 1,\;\;\mathrm{if}\;P_{i}^k/P_{i}^{rate} > \frac{{t + {t_{{\rm{lock}}}}}}{{t + {t_{{\rm{lock}}}} + \Delta t + {t_{{\rm{lock}}}}}}\\
{u_1} = 0.005,\;\; \mathrm{if}\;0.5{\rm{ < }}P_{i}^k/P_{i}^{rate} \le \frac{{t + {t_{{\rm{lock}}}}}}{{t + {t_{{\rm{lock}}}} + \Delta t + {t_{{\rm{lock}}}}}}\\
{u_0} = 0.005,\;\; \mathrm{if}\;0.5\ge P_{i}^k/P_{i}^{rate} \ge \frac{{\Delta t + {t_{{\rm{lock}}}}}}{{{\Delta t}+ {t_{{\rm{lock}}}} + t + {t_{{\rm{lock}}}}}}\\
{u_0} = 1,\;\;\mathrm{if}\; P_{i}^k/P_{i}^{rate} < \frac{{\Delta t + {t_{{\rm{lock}}}}}}{{\Delta t + {t_{{\rm{lock}}}} + t + {t_{{\rm{lock}}}}}}
\end{array} \right.\,.
\end{equation}

By far, the transition probabilities ${{u_0}}$ and ${{u_1}}$ can be determined by Equation \eqref{eq:value_of_u0_and_u1}, with smooth transient performance guaranteed.

\section{Autonomous Control of Air-conditioners}
\label{sec:Physical Model}
This section further applies the proposed SMM to the control of air-conditioners (ACs)..  
\subsection{Physical Model}
The thermodynamic model of an AC system is used to describe the thermal dynamic transition process in a room. The simplified equivalent thermal parameters (ETP) model \cite{He2014Aggregate, song2018thermal} of a cooling AC system describes the dynamic behaviors of indoor temperature ${T_{\rm{a}}}$ with cooling capacity ${Q_{AC}^k}$ as follows:
\begin{equation}
 \label{eq:thermodynamic_model}
 {C_a}\frac{{d{T_a}}}{{dt}} = - \frac{1}{{{R_a}}}\left( {{T_a} - T_o^k} \right) - m\left( t \right){Q_{AC}^k}\,,
\end{equation}
where $T_{\rm{o}}^k$ represents the outdoor temperature at period $k$; ${C_{\rm{a}}}$ and ${{R_{\rm{a}}}}$ denote the equivalent capacitance and resistance, respectively; $t$ is the interval time of period.

The exponential solution of Equation \eqref{eq:thermodynamic_model} is
\begin{equation}
\label{eq:temperature_equation}
T_{\rm{a}}^{k + 1} = T_{\rm{o}}^k + {Q_{AC}^k}{R_a} - (T_{\rm{o}}^k + {Q_{AC}^k}{R_a} - T_{\rm{a}}^k){e^{ - \frac{{dt}}{{{R_{\rm{a}}}{C_{\rm{a}}}}}}}\,.
\end{equation}

Denote $\rm COP$ as the energy efficiency ratio of the ON/OFF controlled AC, then the relation of the dynamic power $P_{{\rm{AC}}}^k$ with $Q_{{\rm{AC}}}^k$ is described as
\begin{equation}
\label{eq:PQrelation_FFA}
Q_{{\rm{AC}}}^k = P_{{\rm{AC}}}^k{\rm{COP}}\,.
\end{equation}
 
Combine \eqref{eq:temperature_equation} and \eqref{eq:PQrelation_FFA}, the physical model of ON/OFF controlled AC is
\begin{equation}
    \label{eq:physical_model}
    T_{\rm{a}}^{k + 1} = T_{\rm{o}}^k +{R_a}  P_{{\rm{AC}}}^k{\rm{COP}}- (T_{\rm{o}}^k + {R_a} P_{{\rm{AC}}}^k{\rm{COP}} - T_{\rm{a}}^k){e^{ - \frac{{dt}}{{{R_{\rm{a}}}{C_{\rm{a}}}}}}}\,.
\end{equation}

According to the SMM, Equation \eqref{eq:expected_value_of_PFFA_equal_response_power} and Equation \eqref{eq:physical_model}, the target power $P_{AC.i}^k$ can be used as a command to control the AC.

\subsection{Control Process Based on Cyber-Physical System (CPS)}

\begin{figure}[H]
 \centering
 \includegraphics[width=12 cm]{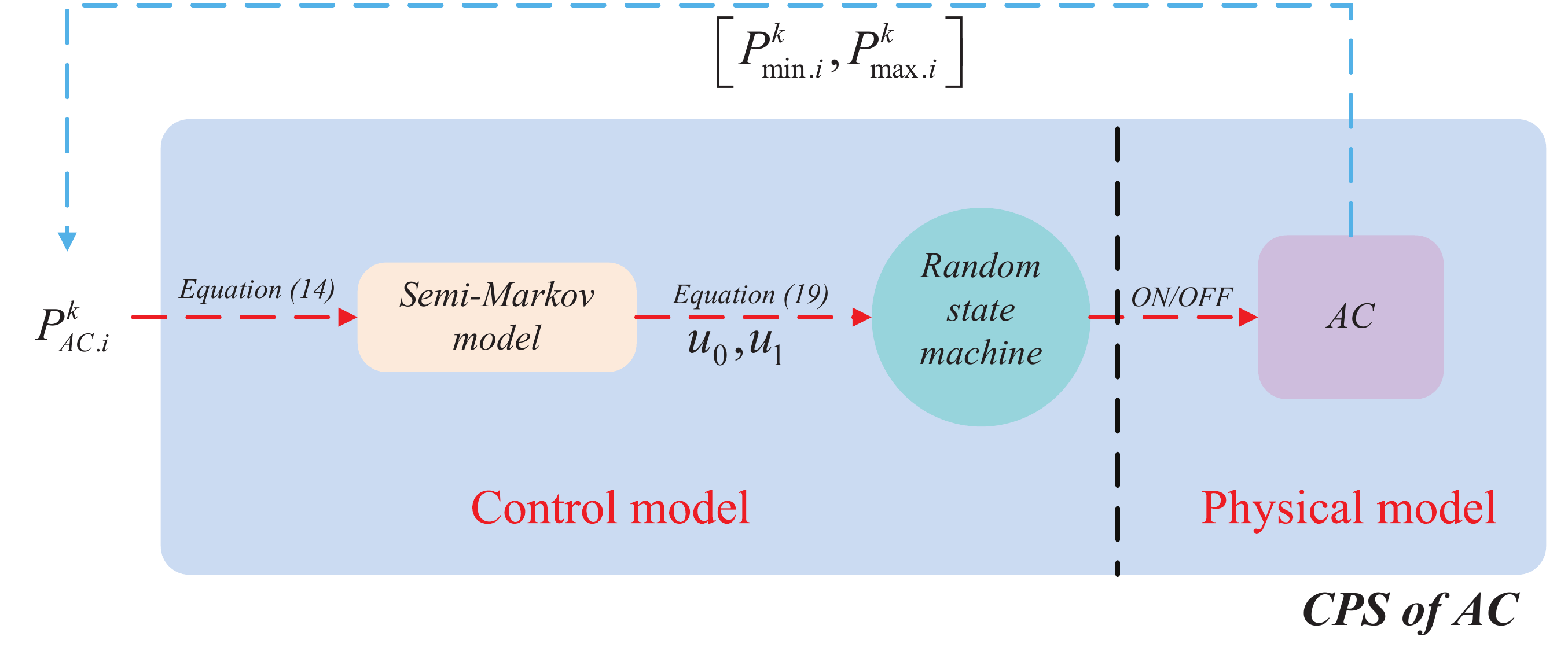}
 \caption{Control process based on CPS.}
 \label{fig:Control_process}
\end{figure}

Combine the SMM model with AC's physical model and the overall control process is demonstrated in Figure \ref{fig:Control_process}. As shown, in control period $k$ the AC $i$ receives a target power $P_{AC.i}^k$ from the aggregator. To track the given target power, $u_0$ and $u_1$ can be obtained consequently from the SMM according to \eqref{eq:value_of_u0_and_u1}. These two transition probabilities are used as input signals to the four-state machine that locally switches the ON/OFF controlled AC within the control period. In this way, the proposed cyber-physical system realizes the autonomous control of each AC by bridging the gap between the continuous target power and the discrete instantaneous power of each AC. Besides, it should be mentioned that both physical constraint and user comfort can be taken into consideration in the proposed system with ease.

Combine \eqref{eq:thermodynamic_model} with \eqref{eq:temperature_equation}, the dynamic cooling capacity $Q_{{\rm{AC}}}^k$ that changes the indoor temperature from ${T_{\rm{a}}^k}$ to ${T_{\rm{a}}^{k+1}}$ within control period $k$ is obtained as
\begin{equation}
\label{eq:Qac_equation}
{Q_{AC}^k} = \frac{{T_{\rm{a}}^{k + 1} - T_{\rm{o}}^k + T_{\rm{o}}^k{e^{ - \frac{{dt}}{{{R_{\rm{a}}}{C_{\rm{a}}}}}}} - T_{\rm{a}}^k{e^{ - \frac{{dt}}{{{R_{\rm{a}}}{C_{\rm{a}}}}}}}}}{{{R_{\rm{a}}} - {R_a}{e^{ - \frac{{dt}}{{{R_{\rm{a}}}{C_{\rm{a}}}}}}}}}\,.
\end{equation}

Therefore, according to \eqref{eq:PQrelation_FFA}and \eqref{eq:Qac_equation},  the dynamic power equation of AC can be calculated as
\begin{equation}
\label{eq:PFFA_equation}
P_{{\rm{AC}}}^k = \frac{{T_{\rm{a}}^{k + 1} - T_{\rm{o}}^k + T_{\rm{o}}^k{e^{ - \frac{{dt}}{{{R_{\rm{a}}}{C_{\rm{a}}}}}}} - T_{\rm{a}}^k{e^{ - \frac{{dt}}{{{R_{\rm{a}}}{C_{\rm{a}}}}}}}}}{{COP\left( {{R_{\rm{a}}} - {R_a}{e^{ - \frac{{dt}}{{{R_{\rm{a}}}{C_{\rm{a}}}}}}}} \right)}}\,.
\end{equation}

According to \eqref{eq:PFFA_equation}, the minimum power $P_{\min .i}^k$ and the maximum power $P_{\max .i}^k$ of AC $i$ during period $k$ can be obtained as
\begin{equation}
\label{eq:PAC_rang}
\left\{ \begin{array}{l}
P_{\min .i}^k = \max \left( {\frac{{T_{a.i}^{\max } - T_o^k + T_o^k{e^{ - \frac{{dt}}{{{R_{a.i}}{C_{a.i}}}}}} - T_{a.i}^k{e^{ - \frac{{dt}}{{{R_{a.i}}{C_{a.i}}}}}}}}{{CO{P_i}\left( {{R_{a.i}} - {R_{a.i}}{e^{ - \frac{{dt}}{{{R_{a.i}}{C_{a.i}}}}}}} \right)}},0} \right)\\
P_{\max .i}^k = \min \left( {\frac{{T_{a.i}^{\min } - T_o^k + T_o^k{e^{ - \frac{{dt}}{{{R_{a.i}}{C_{a.i}}}}}} - T_{a.i}^k{e^{ - \frac{{dt}}{{{R_{a.i}}{C_{a.i}}}}}}}}{{CO{P_i}\left( {{R_{a.i}} - {R_{a.i}}{e^{ - \frac{{dt}}{{{R_{a.i}}{C_{a.i}}}}}}} \right)}},P_{AC.i}^{rate}} \right)
\end{array} \right.\,.
\end{equation}

At the end of every control period, the actual adjustment range $[P_{\min .i}^k,\;P_{\max .i}^k]$of each AC should be uploaded accordingly before the control proceeds.

\section{Case Studies}
\label{sec:case_studies}
\subsection{Parameters Settings}
Table \ref{tab:parameters_of_AC} shows the parameter distributions of a cluster of heterogenous ACs \cite{Yao2019Aggregating}, where $U$ denotes the uniform distribution within the given range. The control interval is $0.5h$ and the execution interval of semi-Markov model is $\Delta t = 2s$.
\begin{table}[htbp]
\centering
\caption{Parameters of AC}
\begin{tabular}{|c|c|c|c|c|c|c|c|}
\hline
Parameter &${R_a}\left( {{\raise0.7ex\hbox{${^\circ C}$} \!\mathord{\left/
 {\vphantom {{^\circ C} {kW}}}\right.\kern-\nulldelimiterspace}
\!\lower0.7ex\hbox{${kW}$}}} \right)$ & ${C_a}\left( {{\raise0.7ex\hbox{${kWh}$} \!\mathord{\left/
 {\vphantom {{kWh} {^\circ C}}}\right.\kern-\nulldelimiterspace}
\!\lower0.7ex\hbox{${^\circ C}$}}} \right)$ & $P_{AC}^{rate}\left( {kW} \right)$ & $COP$ & ${T_{lock}}\left( {\min } \right)$ & ${T_{a.i}^{\min }}\left( {^\circ C} \right)$ & ${T_{a.i}^{\max }}\left( {^\circ C} \right)$ \\\hline
Value & $U\left( {2.5,3.5} \right)$ & $U\left( {1.5,2.5} \right)$ & $U\left( {2.5,3} \right)$ & $U\left( {2.5,3} \right)$ & 3 & $23$ & $27$ \\\hline
\end{tabular}
\label{tab:parameters_of_AC}
\end{table}

\subsection{Distribution of states}

To simulate the stability performance of the proposed model in coordination large-scale TCLs, two scenarios are considered with a cluster of 1,000 and 10,000 ACs respectively. Both parameters and initial states of all ACs are set to the same deliberately. The transition probabilities of this simulation is ${u_0} = 0.0075$ and ${u_1} = 0.0012$ and corresponding theoretical probabilities of stability are ${P_1} = 0.119$, ${P_2} = 0.719$, ${P_3} = 0.081$ and ${P_4} = 0.081$ according to \eqref{eq:probability_of_stability}.

\begin{figure}[H]
 \centering
 \includegraphics[width=12 cm]{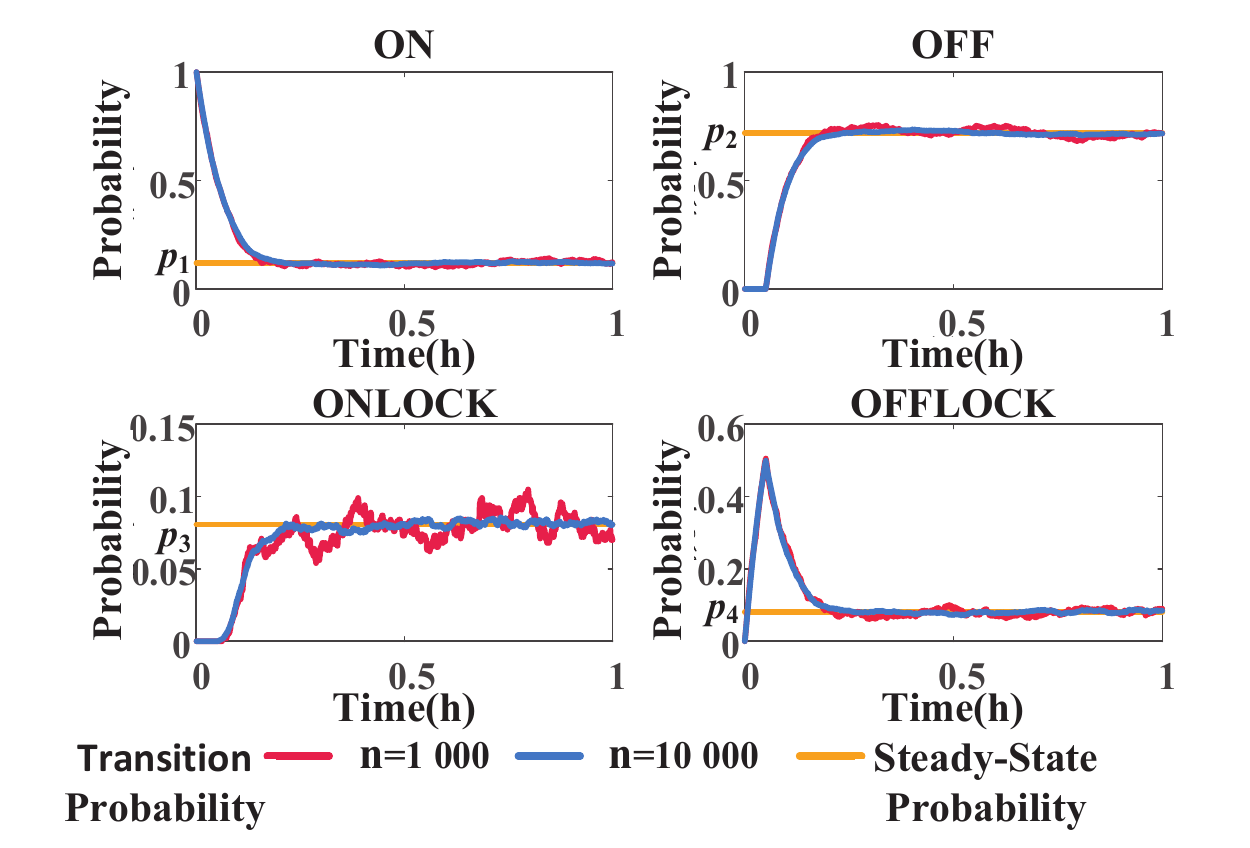}
 \caption{Transient probability of different states.}
 \label{fig:probability_of_stability}
\end{figure}

Each state's actual proportions within the cluster, i.e., the transient probabilities, are plotted in Figure \ref{fig:probability_of_stability}. It can be seen that the transient probability would converge to the theoretical value within an epoch less than $0.5h$. Besides, with the given probability control method, the diversity of operation state within the cluster is guaranteed and the power oscillations in the transient phase are rather tolerable. Furthermore, the cluster with a larger scale has milder fluctuation for all states during the transient phase.

\subsection{Response Performance}
To measure the response performance of the semi-Markov model, this case simulates the coordination of 1000 ACs where the control signals are generated randomly within the limit of operation and user comfort. The simulation lasts for $24h$.

Therefore the random response power $P_{rand.i}^k$ is
\begin{equation}
\label{eq:random_PFFA}
P_{rand.i}^k = rand\left( {P_{\min .i}^k,P_{\max .i}^k} \right)\,.
\end{equation}

\subsubsection{User comfort}
To reflect each user's comfort, the state of AC ($SOA$) is defined as
\begin{equation}
\label{eq:SOA}
SOA_i^k = \frac{{T_{a.i}^k - T_{a.i}^{\min }}}{{T_{a.i}^{\max } - T_{a.i}^{\min }}}\,,
\end{equation}
where $\left[ {T_{a.i}^{\min },T_{a.i}^{\max }} \right]$ is a user specified range of indoor temperature and $SOA_i^k \in [0,1]$.

\begin{figure}[H]
 \centering
 \includegraphics[width=12 cm]{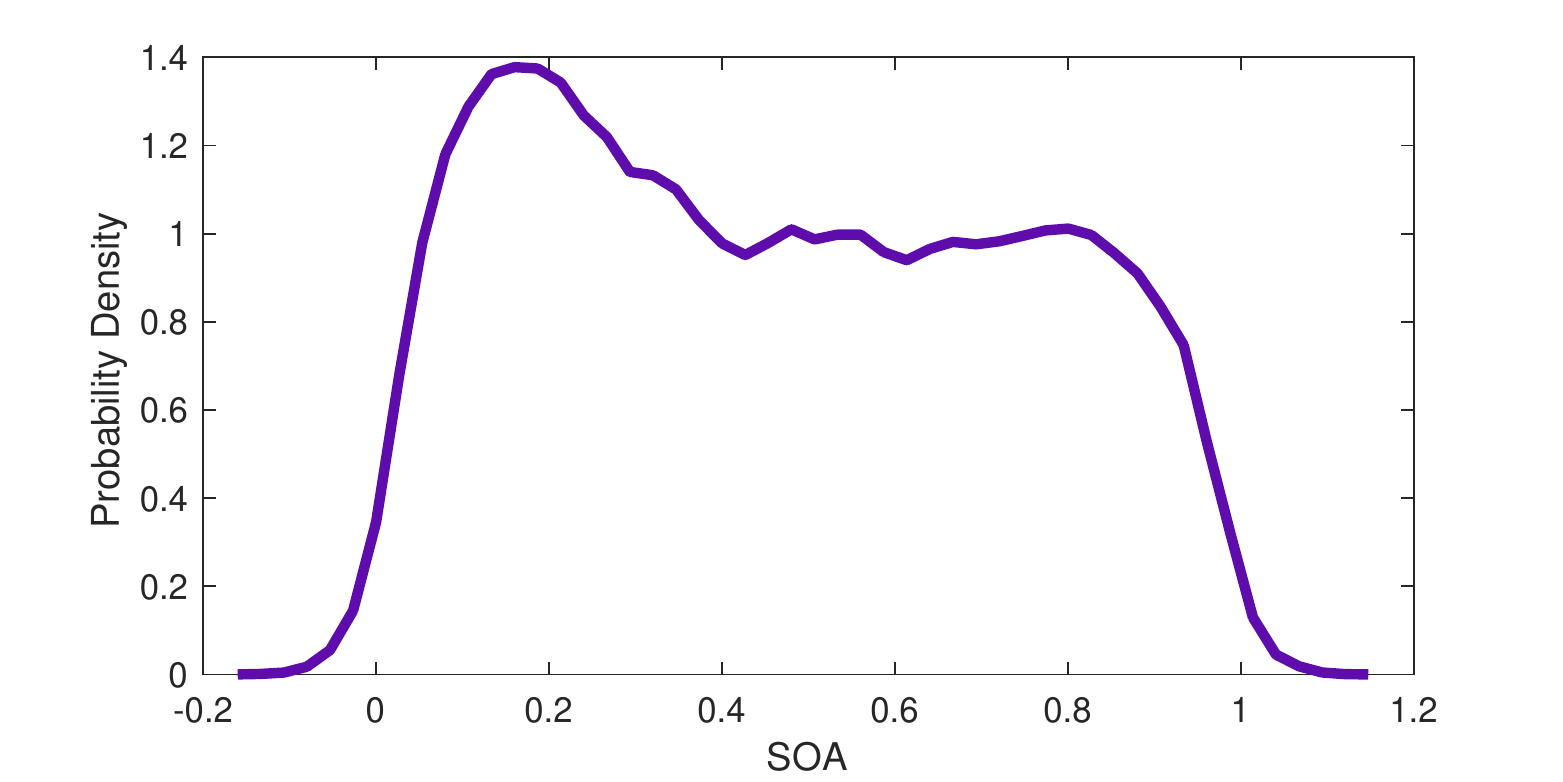}
 \caption{The probability density of SOA.}
 \label{fig:SOA}
\end{figure}

Figure \ref{fig:SOA} plots the probability density of the SOA within the cluster. As shown, the $SOA$ stays within $\left[ {0,1} \right]$ for most of the time while it exceeds the comfort limit slightly due to the inherent uncertainty of the probability control.

\subsubsection{Tracking Performance}
To analyze the performance of the cluster as a whole in tracking the target power, the error of power response is defiend as
\begin{equation}
\label{eq:Er}
{Error}^k{\rm{ = }}\frac{{\sum\limits_{i = 1}^N {\tilde P_{AC.i}^k} - \sum\limits_{i = 1}^N {P_{AC.i}^k} }}{{\sum\limits_{i = 1}^N {P_{AC.i}^{rate}} }}\,,
\end{equation}
where $N$ is the number of AC in the cluster and $P_{AC.i}^k$ is the target power of AC $i$ during $k$ period.
\begin{figure}[H]
 \centering
 \includegraphics[width=12 cm]{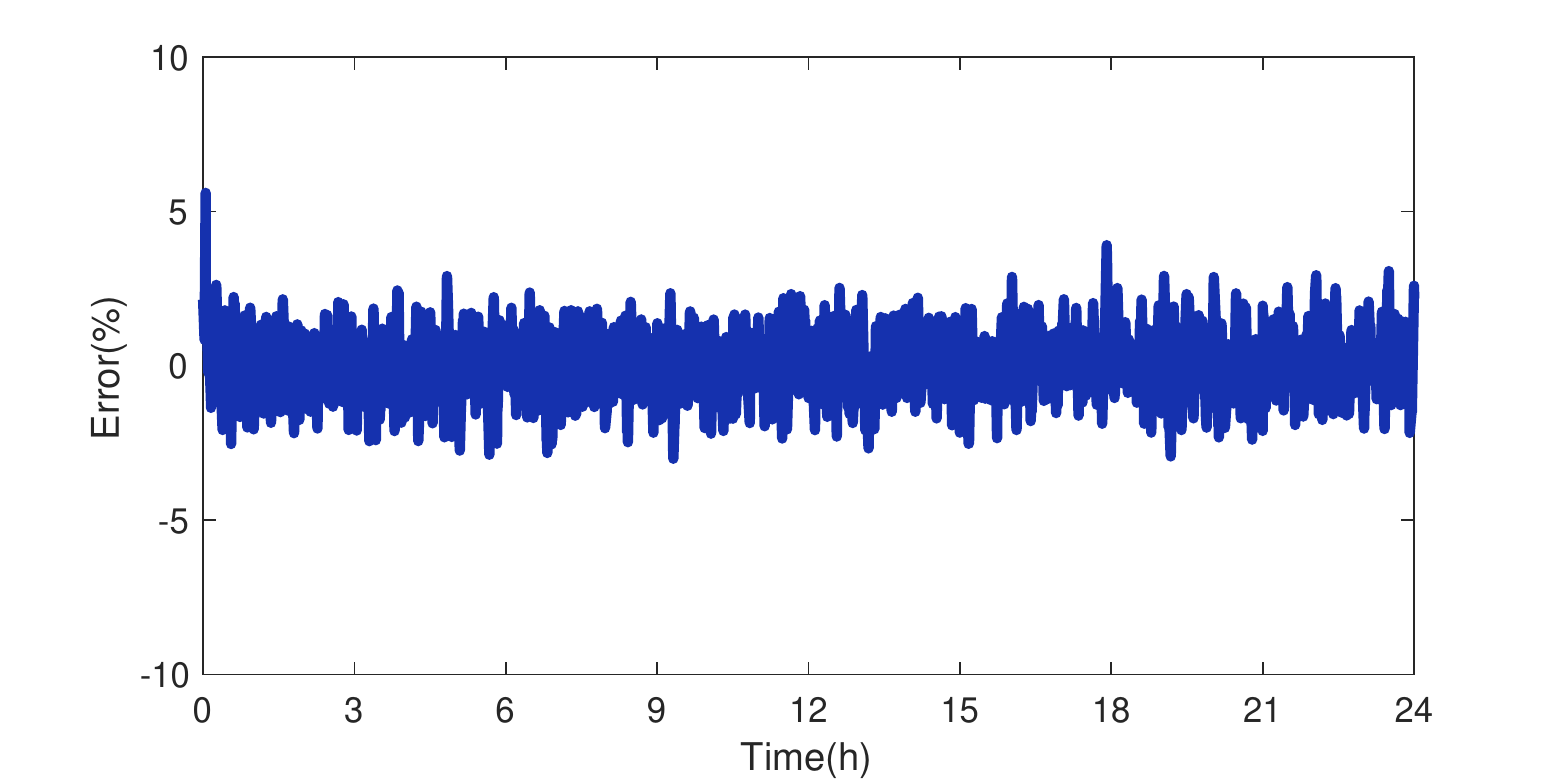}
 \caption{Performance of power response.}
 \label{fig:Power_response}
\end{figure}
Figure \ref{fig:Power_response} shows the performance of power response. As shown, the power deviation is within an acceptable range and the cluster can track the target well with the proposed method.

\section{Conclusions}
\label{sec:conclusions}
To coordinate heterogenous TCLs, frequently switches are typically required to guarantee individual user comfort while tracking a given target power. Considering existing computation and communication infrastructures, this paper proposes a semi-Markov model together with the cyber-physical system for autonomous control of TCL. The proposed model is especially efficient for coordinating large-scale TCLs with the following advantages: (1) The lock time is incorporated in the SMM and the cluster's diversity of switch state is guaranteed; (2) The proposed CPS implements a complete description of two different types of processes, continuous physical dynamics and discrete control logic. (3) Both user comfort and privacy are protected.

\vspace{6pt}




\funding{This work was supported by National Key R\&D Program of China (Basic Research Class 2017YFB0903000).}


\conflictsofinterest{The authors declare no conflict of interest}

\reftitle{References}
\bibliography{ref}

\begin{thebibliography}{-------}
\providecommand{\natexlab}[1]{#1}

\bibitem[Lu and Chassin(2004)]{lu2004state}
Lu, N.; Chassin, D.P.
\newblock A state-queueing model of thermostatically controlled appliances.
\newblock {\em IEEE Transactions on Power Systems} {\bf 2004}, {\em
  19},~1666--1673.

\bibitem[Hao \em{et~al.}(2014)Hao, Sanandaji, Poolla, and
  Vincent]{hao2014aggregate}
Hao, H.; Sanandaji, B.M.; Poolla, K.; Vincent, T.L.
\newblock Aggregate flexibility of thermostatically controlled loads.
\newblock {\em IEEE Transactions on Power Systems} {\bf 2014}, {\em
  30},~189--198.

\bibitem[Chengshan \em{et~al.}(2012)Chengshan, Mengxuan, NING,
  et~al.]{WANGChengshan2012}
Chengshan, W.; Mengxuan, L.; NING, L.; others.
\newblock A Tie-line Power Smoothing Method for Microgrid Using Residential
  Thermostatically-controlled loads.
\newblock {\em Proceedings of the CSEE} {\bf 2012}, pp. 36--43.

\bibitem[Dan \em{et~al.}(2014)Dan, Menghua, and Hongjie]{WANGDan2014}
Dan, W.; Menghua, F.; Hongjie, J.
\newblock User Comfort Constrain Demand Response for Resident
  Thermostically-controlled Loads and Efficient Power Plant Modeling.
\newblock {\em Proceedings of the CSEE} {\bf 2014}, {\em 34},~2071--2077.

\bibitem[Paccagnan \em{et~al.}(2015)Paccagnan, Kamgarpour, and
  Lygeros]{paccagnan2015range}
Paccagnan, D.; Kamgarpour, M.; Lygeros, J.
\newblock On the range of feasible power trajectories for a population of
  thermostatically controlled loads.
\newblock  2015 54th IEEE conference on decision and control (CDC). IEEE,
  2015, pp. 5883--5888.

\bibitem[Mathieu \em{et~al.}(2012)Mathieu, Koch, and
  Callaway]{mathieu2012state}
Mathieu, J.L.; Koch, S.; Callaway, D.S.
\newblock State estimation and control of electric loads to manage real-time
  energy imbalance.
\newblock {\em IEEE Transactions on Power Systems} {\bf 2012}, {\em
  28},~430--440.

\bibitem[Zhang \em{et~al.}(2012)Zhang, Kalsi, Fuller, Elizondo, and
  Chassin]{zhang2012aggregate}
Zhang, W.; Kalsi, K.; Fuller, J.; Elizondo, M.; Chassin, D.
\newblock Aggregate model for heterogeneous thermostatically controlled loads
  with demand response.
\newblock  2012 IEEE Power and Energy Society General Meeting. IEEE,  2012, pp.
  1--8.

\bibitem[Totu \em{et~al.}(2016)Totu, Wisniewski, and Leth]{totu2016demand}
Totu, L.C.; Wisniewski, R.; Leth, J.
\newblock Demand response of a TCL population using switching-rate actuation.
\newblock {\em IEEE Transactions on Control Systems Technology} {\bf 2016},
  {\em 25},~1537--1551.

\bibitem[Moya \em{et~al.}(2014)Moya, Zhang, Lian, and
  Kalsi]{moya2014hierarchical}
Moya, C.; Zhang, W.; Lian, J.; Kalsi, K.
\newblock A hierarchical framework for demand-side frequency control.
\newblock  2014 American Control Conference. IEEE,  2014, pp. 52--57.

\bibitem[Williams \em{et~al.}(2016)Williams, Kalsi, Elizondo, Marinovici, and
  Pratt]{williams2016control}
Williams, T.; Kalsi, K.; Elizondo, M.; Marinovici, L.; Pratt, R.
\newblock Control and coordination of frequency responsive residential water
  heaters.
\newblock  2016 IEEE Power and Energy Society General Meeting (PESGM). IEEE,
  2016, pp. 1--5.

\bibitem[Kalsi \em{et~al.}(2015)Kalsi, Hansen, Fuller, Marinovici, Elizondo,
  Williams, Lian, and Sun]{kalsi2015loads}
Kalsi, K.; Hansen, J.; Fuller, J.C.; Marinovici, L.D.; Elizondo, M.A.;
  Williams, T.L.; Lian, J.; Sun, Y.
\newblock Loads as a resource: Frequency responsive demand.
\newblock Technical report, Pacific Northwest National Lab.(PNNL), Richland, WA
  (United States),  2015.

\bibitem[Harlamov(2008)]{harlamov2008continuous}
Harlamov, B.
\newblock {\em Continuous semi-Markov processes}; Wiley Online Library,  2008.

\bibitem[Zhang \em{et~al.}(2013)Zhang, Lian, Chang, and
  Kalsi]{zhang2013aggregated}
Zhang, W.; Lian, J.; Chang, C.Y.; Kalsi, K.
\newblock Aggregated modeling and control of air conditioning loads for demand
  response.
\newblock {\em IEEE transactions on power systems} {\bf 2013}, {\em
  28},~4655--4664.

\bibitem[He \em{et~al.}(2014)He, Sanandaji, Poolla, and
  Vincent]{He2014Aggregate}
He, H.; Sanandaji, B.M.; Poolla, K.; Vincent, T.L.
\newblock Aggregate Flexibility of Thermostatically Controlled Loads.
\newblock {\em IEEE Transactions on Power Systems} {\bf 2014}, {\em
  30},~189--198.

\bibitem[Song \em{et~al.}(2018)Song, Gao, Yan, and Yang]{song2018thermal}
Song, M.; Gao, C.; Yan, H.; Yang, J.
\newblock Thermal battery modeling of inverter air conditioning for demand
  response.
\newblock {\em IEEE Transactions on Smart Grid} {\bf 2018}, {\em
  9},~5522--5534.

\bibitem[Yao \em{et~al.}(2019)Yao, Zhang, and Chen]{Yao2019Aggregating}
Yao, Y.; Zhang, P.; Chen, S.
\newblock Aggregating Large-Scale Generalized Energy Storages to Participate in
  Energy Market and Regulation Market.
\newblock {\em Energies} {\bf 2019}, {\em 12},~1024.

\end{thebibliography}
\end{document}